\documentclass[12pt]{article}

\usepackage{amsmath}
\usepackage{amsthm}
\usepackage{amsfonts}
\usepackage{latexsym}
\usepackage{epsfig}
\usepackage{fancyhdr}
\usepackage{plain}
\usepackage{eurosym}
\usepackage{mathrsfs}

\newcommand{\beq} {\begin{equation}}
\newcommand{\eeq} {\end{equation}}

\begin{document}
\thispagestyle{empty}
\begin{center}
\textbf{\Large The 'Core' of Symmetric Homogeneous Polynomial Inequalities of Degree Four of Three Real Variables} 
\end{center}
\vspace{1cm}
\noindent \textbf{Mariyan Milev $^{1}$, Nedelcho Milev $^{2}$}\\

\noindent $^1$  Department of Mathematics and Physics, UFT-Plovdiv, bul. Maritza 26, 4000 Plovdiv, Bulgaria, 
tel. +359 32 603744 email: marianmilev2002@gmail.com \\

\noindent $^2$ Department of Mathematical Analysis, Plovdiv University 'Paisii Hilendarski', Tzar Asen 24, 4000 Plovdiv, Bulgaria,
\noindent phone: +359-32 519542, email: milev0@abv.bg \\

%\ \hrulefill \\\
\vspace{0.3cm}
%\hrulefill \\
\noindent\textbf{Abstract.} In this paper we explore inequalities between symmetric homogeneous polynomials of degree four of three real variables and three nonnegative real variables. The main theorems describe the cases in which the smallest possible coefficient is not expressed by the other coefficients. The problem is resolved by introducing a parametric representation.

\vspace{0.3cm}

\vspace{0.3cm}
%\noindent \textbf{MSC}: 65C20

%\newpage
\begin{center}
\Large{\textbf{Introduction}}
\end{center}

In recent years, inequalities between symmetric homogeneous polynomials are studied in several articles \cite{Ando}-\cite{Zhou}. We will explore the inequality
%In this paper we will find necessary and sufficient conditions for validity of cyclic homogeneous polynomials of degree four 
\beq x^4+y^4+z^4 +a(xy^3+yz^3+zx^3+x^3y+y^3z+z^3x)+b(x^2y^2+y^2z^2+z^2x^2)+\nonumber \eeq 
\beq +c \: xyz(x+y+z) \geq 0, \nonumber \eeq  
where $a$, $b$, $c$ are real constants, and $x$, $y$, $z$ are real variables or nonnegative real variables. The proof of the main theorem 3 in this paper is structured by setting of auxiliary functions (2), lemma 4 and the representation (7). For brevity we will denote the main symmetric homogeneous polynomials of degree four 
 \beq \ \ \ w_4=x^4+y^4+z^4, \ \ \ w_3=x^3y+y^3z+z^3x+xy^3+yz^3+zx^3, \ \ w_2=x^2y^2+y^2z^2+z^2x^2 \nonumber \eeq
%\beq  \ \ w_1=xyz(x+y+z) \ \ \ \text{and}. \nonumber \eeq
$w_1=xyz(x+y+z)$ and the inequality is written in the form
\beq f(x,y,z)=w_4+aw_3+bw_2+cw_1\geq 0.  \ \ \ \ \ \ \ \ \ \ \ \ \ \ \ \ \ \ \ \ \ \ \ \ \ \ \ \ \ \ \ \ \ \ \ \ \ \ \ \ \ \ \ \ \ \nonumber\eeq

In Section I we find $b_{\: min}$ such that the inequality $f(x,y,z)=w_4+aw_3+bw_2+cw_1\geq 0$ holds for arbitrary $x,y,z\in R$ if and only if when $b\geq b_{\: min}$.

In the particular case $a=-\displaystyle\frac{1}{2}$ the upper inequality holds if and only if when \\ $b_{\: min}=\max\Big(-c\:;\displaystyle\frac{c}{2}-\displaystyle\frac{9}{8}\Big)$ (Theorem 2).

When $c \leq -a^2-2a$ (Theorem 1) \ $b_{\: min}=-2a-c-1$. This statement follows also from \cite{Cirtoaje}, Theorem 2.1.

When $c>-a^2-2a$ the situation radically changes. In some particular cases $b_{\: min}$ could be expressed in radicals, for example when $a=2$ and $c=4>-a^2-2a$ we find \\ $b_{\: min}=\displaystyle\frac{\Big(5\sqrt{5}-7\Big)}{2}$. When $b=b_{\: min}$ the equality holds when $x=y=1$ and $z=-\displaystyle\frac{1+\sqrt{5}}{2}$ and of course for any permutation.
When $a=-1$ and $c=2>-a^2-2a$ we find $b=\displaystyle\frac{2+5\sqrt[3]{100}-10\sqrt[3]{10}}{12}=0.305299773...$. The equality holds when $x=y=1$ and $z=\displaystyle\frac{\sqrt[3]{100}-2\sqrt[3]{10}-2}{6}=-0.27788...$.

We overcome this problem by introducing the rational parameterization \eqref{two} (see Theorem 3).

We note that this parametrization is in line with \cite{Cirtoaje2}, Theorem 2.2: 'The inequality $f(x,y,z)=w_4+aw_3+bw_2+cw_1\geq 0$ holds if and only if when $f(x,1,1)\geq 0$ for every real $x$.' \\

In Section II we will explore the inequality
\beq f(x,y,z)=w_4+aw_3+bw_2+cw_1\geq 0\: ,\ \ \ \ \ \ \ \ \ \ \ \ \ \ \ \ \ \ \ \ \ \ \ \ \ \ \ \ \ \ \ \ \ \ \ \ \ \ \ \ \ \ \ \ \ \nonumber\eeq
where $x$, $y$, $z$ are nonnegative real variables. %, and $a$, $b$, $c$ are real constants. \\ 

First we find the necessarily conditions for the coefficient $b$ (Lemma 5) and then we search $c_{min}$ such that the previous inequality holds if and only if when $c\geq c_{min}$. We note that $w_1$ is the smallest symmetric polynomial $2w_4 \geq w_3 \geq 2w_2 \geq 2w_1 \geq 0$.

When $a\geq -1$ $\Rightarrow$ $b\geq -2(a+1)$ and $c_{min}=-2a-b-1$. This statement follows also from \cite{Cirtoaje}, Theorem 2.6. According to this theorem the inequality $f(x,y,z)\geq 0$ holds also when $a<-1$, but when $c>-2a-b-1$ the condition $b\geq a^2-1$ is not necessary. When $a<-1$ and $b<a^2-1$, only in some particular cases $c_\emph{min}$ could be expressed in radicals. For example, when $a=-4$ and $b=14$ ($-2(a+1)<b<a^2-1$)
$c_{min}=4\Big[-9+8\sqrt{2}-2\Big(3-\sqrt{2}\Big)\sqrt{3-\sqrt{2}}\Big]=-6.72076...$ When $c=c_{min}$ the equality holds when $x=y=1$ and $z=\sqrt{2}+\sqrt{6-2\sqrt{2}}=3.1951...$ and of course for any permutation. Another example, when $a=-6$, $b=31$ \\ $c_{min}=-58+20\sqrt{15}+\Big(5-3\sqrt{15}\Big)\sqrt{40-2\sqrt{15}}=-18.131094...$
$\Big(\displaystyle\frac{a^2}{4}+2<b<a^2-1 \Big)$. The equality holds when $x=y=1$ and $z=\displaystyle \frac{1+\sqrt{15}+\sqrt{40-2\sqrt{15}}}{2}=5.276...$

We again solve the problem by using the rational parametrization \eqref{two} (see Theorem 7). 

We note that this do not contradict to \cite{Ando2} and \cite{Cirtoaje3}, Theorem 2.1: 'The inequality $f(x,y,z)=w_4+aw_3+bw_2+cw_1\geq 0$ holds if and only if when $f(x,1,0)\geq 0$ and $f(x,1,1)\geq 0$ for every $x\geq 0$.'  

\textbf{Remark.} All identities are verified via the Maplesoft platform.
\begin{center}
\Large{\textbf{Main Results}}
\end{center}

\textbf{Section I.} We will explore the inequality
\beq f(x,y,z)=w_4+aw_3+bw_2+cw_1\geq 0\: ,\label{one} \eeq 
where $a$, $b$, $c$ are real constants and $x$, $y$, $z$ are \textbf{real} variables.  

\newpage \textbf{Theorem 1.} One necessary condition to be true the inequality \eqref{one} is $b \geq -2a-c-1$. When $c \leq -a^2-2a$ this condition is also a sufficient condition (This statement follows also from \cite{Cirtoaje}, Theorem 2.1). 

\textbf{Theorem 2.} When $a=-\displaystyle\frac{1}{2}$ the inequality \eqref{one} holds if and only if when \\
$b \geq \max\Big(-c, \; \displaystyle\frac{c}{2}-\displaystyle\frac{9}{8} \;\Big).$  

When $a \neq-\displaystyle\frac{1}{2}$ and $c \geq -a^2-2a$ we introduce the rational functions:
\beq c(t) =\frac{2t^4+2(a+1)t^3+(2a+4)t^2+(3a+4)t-a}{2t+1} \: ,\ \ \ \ \ \ \ \ \ \ \ \ \ \ \ \ \ \ \ \ \ \ \ \ \ \  \ \ \ \ \ \ \ \ \ \ \ \ \ \ \ \ \ \ \ \  \ \ \ \ \ \ \ \ \ \ \label{two} \eeq 
\beq b(t) =\frac{t^4+(a+4)t^3+(5a+4)t^2+4(a+1)t+2a+2}{-(2t+1)} \ \ \ \ \ \ \ \ \ \ \ \ \  \ \ \ \ \ \ \ \ \ \ \ \ \ \ \ \ \ \ \ \  \ \ \ \ \ \ \ \ \ \ \nonumber \eeq 
and we will explore the different cases about the parameter $a$ when the rational function $c(t)$ always takes values in the interval $[-a^2-2a\: ; \:+\infty)$. \\

\noindent When $-\displaystyle\frac{1}{2}<a \leq 4$ and $t \in [-a-1;-\displaystyle\frac{1}{2})$ the continuous function $c(t)$ is monotone increasing and $-a^2-2a=c(-a-1)\leq c(t)<c(-\displaystyle\frac{1}{2}-0)= +\infty$. \\

\noindent When $-2 \leq a <-\displaystyle\frac{1}{2}$ and $t \in (-\displaystyle\frac{1}{2} ; -a-1]$ the continuous function $c(t)$ is monotone decreasing and 
$+\infty=c(-\displaystyle\frac{1}{2}+0)>c(t) \geq c(-a-1) =-a^2-2a$. \\

\noindent Let us set
\beq t_1=\frac{-a-\sqrt{(a+2)(a-4)}}{2} \ \ \ \ \ \text{and} \ \ \ \  t_2=\frac{-a+\sqrt{(a+2)(a-4)}}{2}\:. \label{three} \eeq
When $a>4$ we obtain two intervals $t\in[-a-1\: ; \: t_1]\cup [\: t_2 \:; \: -\displaystyle\frac{1}{2}\Big)$, \\ $-a-1<t_1<-2<t_2<-\displaystyle\frac{1}{2}$. The function $c(t)$ is monotone increasing and 
\beq -a^2-2a=c(-a-1)<c(t_1)=c(t_2)=\frac{a^2}{2}+2<c(-\frac{1}{2}-0) =+\infty \:. \ \ \ \ \ \ \ \ \ \ \ \ \ \ \ \ \ \ \ \ \ \ \ \ \ \ \nonumber\eeq
When $a<-2$ we obtain again two intervals: $t\in (-\displaystyle\frac{1}{2} \: ; \: t_1]\cup[-a-1 \:; \: t_2]$. The function $c(t)$ is monotone decreasing in $(-\displaystyle\frac{1}{2} \: ; \: t_1]$ and monotone increasing in $[-a-1 \:;\:t_2]$:
\beq +\infty=c(-\displaystyle\frac{1}{2}+0)>c(t_1)=c(t_2)=\frac{a^2}{2}+2>c(-a-1)=-a^2-2a \:. \ \ \ \ \ \ \ \ \ \ \ \ \ \ \ \ \ \ \ \ \ \ \ \ \ \ \ \ \nonumber \eeq
\textbf{Theorem 3.} When $a\neq -\displaystyle\frac{1}{2}$ and $c\geq -a^2-2a$ in the following four cases: \\
(3.1.) \ \ \ $-\displaystyle\frac{1}{2} <a \leq 4$ and $t\in [-a-1 \: ; -\displaystyle\frac{1}{2})$; \\
(3.2.) \ \ \ $ -2 \leq a <-\displaystyle\frac{1}{2}$ and $t\in (-\displaystyle\frac{1}{2}\: ; -a-1 ]$;\\
(3.3.) \ \ \ \ \ $a>4$ and $t\in[-a-1 \: ;t_1]\cup [t_2 \:; -\displaystyle\frac{1}{2}\Big)$;\\
(3.4.) \ \ \ \ \ $a<-2$ and $t\in \Big(-\displaystyle\frac{1}{2} \: ; t_1]\cup[-a-1\:;t_2]$ \\
the inequality $ f(x,y,z)=w_4+aw_3+bw_2+cw_1\geq 0$ holds if and only if when $b\geq b(t)$. When $b = b(t)$ we have $f(t,1,1)=f(1,t,1)=f(1,1,t)=0$. \\

For the proof of this theorem 3 and theorem 7 below we will use the following known statement. 

\textbf{Lemma 4.} For arbitrary real numbers $X$, $Y$ and $Z$ the inequality
\beq g(X, Y, Z)=X^2+Y^2+Z^2+k(XY+YZ+ZX) \geq 0\label{four} \eeq
holds if and only if when $-1\leq k \leq 2$. \\

\textbf{Section II.} We will explore the inequality
\beq f(x,y,z)=w_4+aw_3+bw_2+cw_1\geq 0\: ,  \ \ \ \ \ \ \ \ \  \ \ \ \ \ \ \ \ \  \ \ \ \ \ \ \ \ \  \ \ \ \ \ \ \ \ \  \ \ \ \ \ \ \ \ \  \ \ \ \ \ \ \ \ \  \label{five} \eeq 
where $a$, $b$, $c$ are real constants and $x$, $y$, $z$ are \textbf{nonnegative} real variables. \\ 

\textbf{Lemma 5.} To be true the inequality \eqref{five} it is necessary:
\begin{displaymath} b \geq \ \left\{
\begin{array}{ll} 
\displaystyle\frac{1}{4}\:a^2+2 & \textrm{when \ $a \leq -4$};\\
-2(a+1) & \textrm{when \ $a \geq -4$} \ \
 \end{array}  \textrm{and} \ \ c\geq -2a-b-1.  \ \ \ \ \ \ \ \ \ \  \ \ \ \ \ \ \ \ \ \ \right.  \ \ \ \ \ \ \ \ \ \  \ \ \ \ \ \ \ \ \ \
\end{displaymath}

\textbf{Theorem 6.} (\cite{Cirtoaje}, Theorem 2.6). When $a\geq -1$ the inequality \eqref{five} holds if and only if when $b \geq -2(a+1)$ and $c\geq -2a-b-1$. 

\textbf{Remark.} According to (\cite{Cirtoaje}, Theorem 2.6) the inequality \eqref{five} is true also when $a<-1$, $b\geq a^2-1$ and $c\geq -2a-b-1$, but when $c> -2a-b-1$
 the condition $b\geq a^2-1$ is not necessary.  

When $a<-1$ and $b\leq a^2-1$ we will describe the different cases about the parameter $a$ again with the rational function $b(t)$ and $c(t)$, defined by \eqref{two}.

When $-2\leq a < -1$ and $t \in [0; -a-1]$ the rational function $b(t)$ is monotone increasing and takes values in the intervala $[-2(a+1) \: ; a^2-1]$, i.e. \\ $-2(a+1)=b(0)\leq b(t) \leq b(-a-1)=a^2-1$. \\

When $-4<a<-2$ there are two intervals for $t$.

When $t\in [0\:;t_1]$ the function $b(t)$ is monotone increasing and takes values in the interval $[-2(a+1)\: ; \displaystyle\frac{a^2}{4}+2]$, i.e. $-2(a+1)=b(0)\leq b(t) \leq b(t_1)=\displaystyle\frac{a^2}{4}+2$, and when $t\in [-a-1\: ;t_2]$ the function $b(t)$ is monotone decreasing and takes values in the interval $[\displaystyle\frac{a^2}{4}+2 \:; a^2-1]$, \: i.e. \: $\displaystyle\frac{a^2}{4}+2=b(t_2)\leq b(t) \leq b(-a-1)=a^2-1$, where $t_1$ and $t_2$ are defined in formulas \eqref{three}. \\

When $a \leq -4$ and $t\in [-a-1;t_2]$ the function $b(t)$ is monotone decreasing and takes values in the interval $[\displaystyle\frac{a^2}{4}+2\: ; a^2-1]$, i.e. $\displaystyle\frac{a^2}{4}+2=b(t_2)\leq b(t) \leq b(-a-1)=a^2-1$. \\

\textbf{Theorem 7.} In the three cases: \\
(7.1.) \ \ \  $-2\leq a <-1$ and $t \in [0\:; -a-1]$ ; \\
(7.2.) \ \ \ $-4<a<-2$ and $t\in [0\:;t_1]\cup [-a-1\:;t_2]$;\\
(7.3.) \ \ \ \ \ $a \leq -4$ and $t\in [-a-1\:;t_2]$\\
the inequality $ f(x,y,z)=w_4+aw_3+b(t)w_2+cw_1\geq 0$ holds if and only if when $c\geq c(t)$. When $c = c(t)$ we have $f(t,1,1)=f(1,t,1)=f(1,1,t)=0$. \\
\begin{center}
\Large{\textbf{Proves of Theorems}}
\end{center}

\textbf{Section I} (real variables) \\

\textbf{Remark.} In the proves we have in mind that the following basic symmetric homogeneous polinomial inequalites of degree four of three real variables are true:
\beq 2w_4 \geq w_3\:, \ \ w_4 \geq w_2\geq 0\:, \ \ w_2 \geq w_1 \ \ \text{and} \ \ w_2+2w_1 \geq 0. \ \ \ \ \ \ \ \ \ \ \ \ \ \ \ \ \ \ \ \ \ \ \ \ \ \ \ \ \ \ \ \ \ \ \ \ \ \ \ \ \nonumber\eeq
\textbf{Proof of Theorem 1.} 

The necessity follows from $0\leq f(1,1,1)=3+6a+3b+3c$, i.e. $b\geq -2a-c-1$.  

The sufficiency follows from the relation
\beq w_4+aw_3+(a^2-1)w_2-(a^2+2a)w_1=u^2+v^2-uv\geq 0, \ \ \ \ \ \ \ \ \ \ \ \ \ \ \ \ \ \ \ \ \ \ \ \ \ \ \ \ \ \ \ \ \ \label{six}\eeq
where $u=(x-y)(x+y+az)$ and $v=(x-z)(x+z+ay)$.

When $c\leq -a^2-2a$ \ $\Rightarrow$ \ $b\geq-2a-c-1\geq a^2-1$ and $w_4+aw_3+bw_2+cw_1=$
\beq =\big[w_4+aw_3+(a^2-1)w_2-(a^2+2a)w_1\big]+(2a+b+c+1)w_2-(c+a^2+2a)(w_2-w_1)\geq 0. \ \ \ \ \ \ \ \ \ \ \ \ \ \ \ \ \ \ \ \ \ \ \ \ \ \ \ \ \ \ \ \ \ \nonumber\eeq
\textbf{Remark.} When $a=-1$, $c=-a^2-2a=1$ and $b=-2a-c-1=0$ we obtain the famous inequality of I. Schur \cite{Schur},  \ $w_4+w_1\geq w_3$. \\

\textbf{Proof Theorem 2.} \\

When $a=-\displaystyle\frac{1}{2}$ and $c\leq -a^2-2a=\displaystyle\frac{3}{4}$ \ $\Rightarrow$ \ $b\geq -2a-c-1=-c\geq -\displaystyle\frac{3}{4}$ the basic inequality \eqref{six} is
$w_4-\displaystyle\frac{1}{2}w_3-\displaystyle\frac{3}{4}w_2+\displaystyle\frac{3}{4}w_1\geq 0$ and the statement follows from Theorem 1.

From $0\leq f(2,2,-1)=24b-12c+27$ \ $\Rightarrow$ \ $b \geq \displaystyle\frac{c}{2}-\displaystyle\frac{9}{8}$.

When $c\geq \displaystyle\frac{3}{4}$ \ and \: $b \geq \displaystyle\frac{c}{2}-\displaystyle\frac{9}{8}$ we have
\beq w_4-\frac{1}{2}w_3+bw_2+cw_1=\Big(w_4-\frac{1}{2}w_3-\frac{3}{4}w_2+\frac{3}{4}w_1\Big)+\frac{1}{2}\Big(c-\frac{3}{4}\Big)(w_2+2w_1)+\Big(b-\frac{c}{2}+\frac{9}{8}\Big)w_2 \geq 0.\nonumber \eeq
Before the proof of Theorem 3 we will prove Lemma 4. \\

\textbf{Proof of Lemma 4.} \\

\noindent From $g(1,1,1)=3+3k\geq 0$ and from $g(1,-1,0)=2-k\geq 0$ it follows that $-1 \leq k \leq 2$. And $g(X,Y,Z)=\displaystyle\frac{k+1}{3}\:\Big(X+Y+Z\Big)^2+\displaystyle\frac{2-k}{3}\:\Big(X^2+Y^2+Z^2-XY-YZ-ZX\Big)\geq 0$. \\

\textbf{Proof of Theorem 3.} Necessity. 
\beq f(x,y,z)=w_4+aw_3+bw_2+c(t)w_1=w_4+aw_3+b(t)w_2+c(t)w_1+\Big(b-b(t)\Big)w_2.\nonumber\eeq
From $0\leq f(t,1,1)=0+\Big(b-b(t)\Big)(2t^2+1)$ follows that $b\geq b(t)$. 

Sufficiency follows from the identity:
\beq w_4+aw_3+b(t)w_2+c(t)w_1=X^2+Y^2+Z^2+k\:\Big (XY+YZ+ZX\Big), \ \text{where}\nonumber \eeq
\beq X=x^2+pxy+pxz+qyz, \ Y=y^2+pyz+pyx+qzx, \ Z=z^2+pzx+pzy+qxy,\label{seven} \eeq
\beq p=-\frac{t^2+t+1}{2t+1}, \ q=\frac{t^2+2t}{2t+1} \ \ \text{and} \ \ k=\frac{2t^2+2t(a+1)+a+2}{t-1}\nonumber \eeq
When $t=1$ we obtain that $p=-1$, $q=1$ and the indentity do not depend on $k$ because $XY+YZ+ZX=0$.

Having in mind Lemma 4 it remains to check that $-1\leq k \leq 2$ for the different cases. \\
(3.1) When $-\displaystyle\frac{1}{2}<a \leq 4$ and $t\in [-a-1\; ; -\displaystyle\frac{1}{2})$ $\Rightarrow$ $2t+1<0$, $t-1<0$, $t+a+1\geq 0$ $\Rightarrow$
\beq k+1=\frac{(2t+1)(t+a+1)}{t-1}\geq 0 \ \ \text{and} \  \ k-2=\frac{(2t+a)^2+(4-a)(a+2)}{2(t-1)}<0.\nonumber \eeq
(3.2) When $-2<a<-1$ and $t\in (-\displaystyle\frac{1}{2}\; ; -a-1]$ $\Rightarrow$ $2t+1>0$, $t-1\leq -a-2 <0$, $t+a+1<0$ \ \ $\Rightarrow$
\beq k+1=\frac{(2t+1)(t+a+1)}{t-1}\geq 0 \ \ \text{and} \ \ k-2=\frac{(2t+a)^2+(4-a)(a+2)}{2(t-1)}<0.\nonumber \eeq
When $a=-2$ and $t\in (-\displaystyle\frac{1}{2}\; ; -a-1)$ \ $\Rightarrow$ \ $k=2t\in(-1;2)$. \\
\textbf{Remark.} The case $t=1$ is possible only if $a=-2$ \ $\Rightarrow$ \ $b(1)=3$, $c(1)=0$ \ and 
\beq 2(w_4-2w_3+3w_2)=(y-z)^4+(x-z)^4+(x-y)^4=2(x^2+y^2+z^2-xy-yz-xz)^2.\nonumber \eeq
(3.3) When $a>4$ and $t\in [-a-1;t_1]\cup[t_2; -\displaystyle\frac{1}{2}\Big)$ $\Rightarrow$ $2t+1<0$, $t-1<0$, $t+a+1\geq 0$ 
\beq \Rightarrow \ \ k+1=\frac{(2t+1)(t+a+1)}{t-1}\geq 0 \ \ \text{and} \ \ k-2=\frac{2(t-t_1)(t-t_2)}{t-1}\leq 0. \ \ \ \ \ \ \ \nonumber \eeq
(3.4) Case: $a<-2$ and $t\in\Big(-\displaystyle\frac{1}{2}\:;t_1]\cup[-a-1\:;t_2]$. \\
When $t\in\Big(-\displaystyle\frac{1}{2}\:;t_1]$ \ $\Rightarrow$ \ $t\leq t_1<1$, \: $2t+1>0$, \: $t+a+1<0$, i.e. 
\beq k+1=\frac{(2t+1)(t+a+1)}{t-1}\geq 0 \ \ \text{and} \ \ k-2=\frac{2(t-t_1)(t-t_2)}{t-1}\leq 0.\nonumber \eeq
When $t\in [-a-1\:;t_2]$ \ $\Rightarrow$ \ $t>1$, \: $t+a+1\geq 0$ \: and again
\beq k+1=\frac{(2t+1)(t+a+1)}{t-1}\geq 0 \ \ \text{and} \ \ k-2=\frac{2(t-t_1)(t-t_2)}{t-1}\leq 0.\nonumber \eeq
\textbf{Remark.} The function $b(t)$ is not monotone, e.g. when \\ 
$a=2$, $t=-1$ \ $\Rightarrow$ \ $c(-1)=8$, $b(-1)=3$, $k=0$, $p=q=1$ \: and 
\beq w_4+2w_3+3w_2+8w_1=(x^2+xy+xz+yz)^2+(y^2+xy+xz+yz)^2+(z^2+xy+xz+yz)^2 \geq 0\nonumber \eeq
while when $t=-3$ \ $\Rightarrow$ \ $c(-3)=-8$, \: $b(-3)=3=b(-1)$, \: $p=\displaystyle\frac{7}{5}$, \: $q=-\displaystyle\frac{3}{5}$ \: and as \: $k=-1$ \: from \eqref{seven} 
\beq w_4+2w_3+3w_2-8w_1=\displaystyle\frac{1}{2} \: \Big[(X-Y)^2+(Y-Z)^2+(Z-X)^2\Big]= \ \ \ \ \ \ \ \ \ \ \ \ \ \ \ \ \ \ \ \ \ \ \ \ \ \ \ \ \ \ \ \ \ \ \ \nonumber \eeq
\beq =\displaystyle\frac{1}{2} \Big[(x-y)^2(x+y+2z)^2+(y-z)^2(y+z+2x)^2+(z-x)^2(z+x+2y)^2 \Big] \geq 0. \ \ \ \ \ \ \ \ \ \ \ \ \ \ \ \ \ \ \ \ \ \ \ \ \ \ \ \ \nonumber \eeq
Analogously when \: $a=\displaystyle\frac{19}{4}$, \: $t=-1$ \: and \: $t=-4$ \: we obtain the inequalities:
\beq w_4+\Big(\frac{19}{4} \Big)\: w_3+\Big(\frac{17}{2} \Big)\: w_2+19w_1 \geq 0 \ \ \ \ \text{and} \ \ \ w_4+\Big(\frac{19}{4} \Big)\:  w_3+\Big(\frac{17}{2} \Big)\: w_2+\Big(\frac{49}{4} \Big)\: w_1 \geq 0. \ \ \ \nonumber \eeq
\textbf{Remark.} The parameter $t$ generates different optimal inequalities. \\
1. When $t=-a-1$ \ we have \
\beq w_4+aw_3+(a^2-1)w_2-(a^2+2a)w_1=  \nonumber \eeq
\beq =\frac{1}{2} \: \Big((x-y)^2(x+y+az)^2+(y-z)^2(y+z+ax)^2+(z-x)^2(z+x+ay)^2\Big).   \nonumber \eeq
2. When \: $t=t_1$ \: or \: $t=t_2$ \: and \: $a\in(-\infty\: ; -2]\cup[4\: ; +\infty)$ \ $\Rightarrow$ \ $b=\displaystyle\frac{a^2}{4}+2$, \: $c=\displaystyle\frac{a^2}{2}+a$, \: $k=2$, \: $2p+q=\displaystyle\frac{a}{2}$ \: and according to \eqref{seven} we have 
\beq w_4+aw_3+\Big(\displaystyle\frac{a^2}{4}+2\Big)w_2-\Big(\displaystyle\frac{a^2}{2}+a\Big)w_1=  \nonumber \eeq
\beq =\Big(X+Y+Z\Big)^2=\Big(x^2+y^2+z^2+\frac{a}{2}\:(xy+yz+zx)\Big)^2. \nonumber \eeq
3. When $t=-2$ \: and \: $a\in[1\: ;4]$ \ we have 
\beq w_4+aw_3+2(a-1)w_2+(5a-8)w_1= \nonumber \eeq
\beq =\Big(x+y+z\Big)^2\Big(x^2+y^2+z^2+(a-2)\:(xy+yz+zx)\Big). \nonumber \eeq
4. When $t=-1$ \: and \: $a\in[0 \:;6]$ \ $\Rightarrow$ \ $p=q=1$, \ $k=-1+\displaystyle\frac{a}{2}$ \ \ and we have
\beq w_4+aw_3+(2a-1)w_2+4aw_1\geq 0.\nonumber \eeq
5. When $t=0$ \: and \: $a\in[-4 \:; -1 ]$ \ $\Rightarrow$ \ $p=-1$, $q=0$, \ $k=-a-2$ \ and we have
\beq w_4+aw_3-2(a+1)w_2-aw_1\geq 0.\nonumber \eeq

\newpage\textbf{Section II.} (nonnegative real variables) 
  
\textbf{Remark.} In the proof we have in mind that the following fundamental symmetric homogeneous polinomial inequalities of degree four for three nonnegative real variables hold:
$2w_4 \geq w_3 \geq 2w_2 \geq w_1 \geq 0$, and the inequality I. Schur \cite{Schur} $w_4+w_1\geq w_3$. 

\textbf{Proof of Lemma 5.} From $f(1,1,1) \geq 0$ \ $\Rightarrow$ \ $c\geq -2a-b-1$. \\
From \ $0\leq f(x,1,0)=x^4+1+a\Big(x^3+x\Big)+bx^2$ \ \ $\Rightarrow$ \ \ $x^2+x^{-2}+a\Big(x+x^{-1}\Big)+b\geq 0$. 
We set \: $u=x+x^{-1}\geq 2$. When $a\geq -4$ \: and \: $u=2$ \ \ $\Rightarrow$ \ \ $b\geq -2(a+1)$. \\
When \: $a\leq -4$ \: and \: $u=-\displaystyle\frac{a}{2}\geq 2$ \ \ $\Rightarrow$ \ \ $b \geq \displaystyle\frac{a^2}{4}+2$. 

\textbf{Proof of Theorem 6.} The necessity follows from Lemma 5. \\
Using the inequalities: \: $w_4+w_1 \geq w_3$, $w_3\geq 2w_2$ \: and \: $w_2\geq w_1\geq 0$ \: we see that the conditions \: $2a+2+b\geq 0$ \: and \: $1+2a+b+c \geq 0$ \: are also sufficient:
\beq f=w_4+aw_3+bw_2+cw_1= (w_4+w_1-w_3)+\nonumber \eeq
\beq+(a+1)(w_3-2w_2)+(2a+2+b)(w_2-w_1)+(2a+b+c+1)w_1\geq 0. \nonumber \eeq

\textbf{Proof of Theorem 7.} Necessity. \\ $f(x,y,z) =w_4+aw_3+b(t)w_2+cw_1=w_4+aw_3+b(t)w_2+c(t)w_1+\Big(c-c(t)\Big)w_1$. \\
From \ $0\leq f(t,1,1)=0+\Big(c-c(t)\Big)t(t+2)$ \: when \: $t>0$ \: follows that \: $c\geq c(t)$.\\
When \ $t=0$ \: then \: $c(0)=-a$, \: $b(0)=-2(a+1)$, \: $a\in (-4\:;-1)$. \\
From \: $0 \leq \displaystyle \lim_{x\rightarrow +0}\frac{f(x,1,1)}{x}=2(c+a)$ \ it follows that \ $c\geq -a = c(0)$. 

Sufficiency follows analogously as in Theorem 3 from the indentity:
\beq w_4+aw_3+b(t)w_2+c(t)w_1=X^2+Y^2+Z^2+k(XY+YZ+ZX),  \nonumber \eeq
defined with formulas \eqref{seven}, where \: $-1\leq k \leq 2$, \: because each of the three cases are subregions of these of Theorem 3: \\
(7.1) $-2 \leq a < -1$ \: and \: $t\in [\: 0\: ; -a-1]$ \: contains itself in\\
(3.2) that is \: $-2 \leq a < -\displaystyle\frac{1}{2}$ \ and \ $t\in \Big(-\displaystyle\frac{1}{2}\: ; -a-1]$; \\
(7.2) $-4 < a < -2$ \: and \: $t\in [0; t_1]\cup [-a-1; t_2]$ \: and \: (7.3) $a\leq -4$ \: and \: $t\in [-a-1; t_2]$ contain itself in (3.4) $a<-2$ \: and \: $t\in (-\displaystyle\frac{1}{2}\:; t_1]\cup [-a-1\:; t_2]$. \\

\end{document}